\pdfoutput=1
\RequirePackage{ifpdf}
\ifpdf 
\documentclass[pdftex]{sigma}
\else
\documentclass{sigma}
\fi

\numberwithin{equation}{section}

\newtheorem{Theorem}{Theorem}[section]
\newtheorem*{Theorem*}{Theorem}

\newtheorem{Lemma}[Theorem]{Lemma}
\newtheorem{Proposition}[Theorem]{Proposition}

\theoremstyle{definition}
\newtheorem{Definition}[Theorem]{Definition}

\begin{document}

\allowdisplaybreaks

\newcommand{\arXivNumber}{2601.19875}

\renewcommand{\PaperNumber}{041}

\FirstPageHeading

\ShortArticleName{Mass, Staticity, and a Riemannian Penrose Inequality for Weighted Manifolds}

\ArticleName{Mass, Staticity, and a Riemannian Penrose Inequality\\ for Weighted Manifolds}

\Author{Stephen MCCORMICK}

\AuthorNameForHeading{S.~McCormick}

\Address{Institutionen f\"or teknikvetenskap och matematik, Lule{\aa} tekniska universitet,\\ 971 87 Lule\aa, Sweden}
\Email{\mail{stephen.mccormick@ltu.se}}
\URLaddress{\url{https://quasilocal.com/}}

\ArticleDates{Received January 29, 2026 in final form April 23, 2026; Published online April 29, 2026}

\Abstract{In this note, we show that the weighted mass of Baldauf and Ozuch (2022) can be derived as a natural geometric mass invariant following Michel (2011), for a certain weighted curvature map. An associated weighted centre of mass definition is also derived from this. The adjoint of the linearisation of this curvature map leads to a notion of weighted static metrics, which are natural candidates for weighted mass minimisers. This weighted curvature quantity is essentially the scalar curvature of a conformally related metric that Law, Lopez and Santiago~(2025) used to considerably simplify the proof of the weighted positive mass theorem. We show an equivalence between static metrics and weighted static metrics via the conformal relationship, from which we show that a uniqueness theorem holds for weighted static manifolds with weighted minimal surface boundaries. Furthermore, we show that weighted manifolds satisfy a Riemannian Penrose inequality whose equality case holds precisely for these unique weighted static metrics.}


\Classification{53C99; 83C99}

\section{Introduction}
The positive mass theorem is a celebrated cornerstone of mathematical general relativity, proven by Schoen and Yau \cite{SchoenYau79} using minimal surface techniques, and by Witten \cite{Witten81} using a spinor argument. Ideas from both proofs, as well as the positive mass theorem itself, have motivated a considerable amount of research in the decades that followed, including many generalisations and extensions. Some generalisations are driven by physical motivations, such as the Penrose inequality or the various inequalities between mass, charge and angular momentum, while other geometric notions of mass have also been studied in relation to scalar curvature problems. One such definition, which we consider here, is the mass of a weighted manifold (manifold with density) introduced by Baldauf and Ozuch \cite{BaldaufOzuch22}.

The weighted manifold $(M^n,g,{\rm e}^{-f}{\rm d}V_g)$ is a Riemannian manifold $(M^n,g)$ equipped with a~smooth function $f$ defining the measure ${\rm e}^{-f}{\rm d}V_g$. The weighted scalar curvature is given by
\[
	R_f = R + 2\Delta f - |\nabla f|^2,
\]
where we use the convention $\Delta=\nabla_i \nabla^i$. Throughout, we assume all manifolds are smooth, connected, and oriented. Baldauf and Ozuch define the \emph{weighted mass} of an asymptotically flat weighted manifold as
\begin{equation*}
	\mathfrak{m}_f(g) = \mathfrak{m}_{\rm ADM}(g)  +  \frac{1}{(n-1) \omega_{n-1}}\lim_{\rho\to\infty}\int_{S_\rho} \nabla_i f \nu^i {\rm e}^{-f} {\rm d}S,
\end{equation*}
where $S_\rho$ are large coordinate spheres in an asymptotically flat end with outward unit normal~$\nu$, $\omega_{n-1}$~is the volume of the unit $(n-1)$-sphere, and
\begin{equation*}
	\mathfrak{m}_{\rm ADM}(g)=\frac{1}{2(n-1) \omega_{n-1}}\lim_{\rho\to\infty}\int_{S_\rho} \bigl( \mathring{\nabla}^jg_{ij}-\mathring{\nabla}_i\operatorname{tr}_{\mathring g}(g)\bigr) \nu^i{\rm d}S
\end{equation*}
is the ADM mass. Here and throughout, we use~$\mathring g$ to denote the Euclidean metric pulled back to the asymptotic end, and $\mathring\nabla$ is the associated flat connection. Note that our definition of the weighted mass differs from that of \cite{BaldaufOzuch22} by a multiplicative constant to align with the standard definition of the ADM mass, so that the resultant Penrose inequality takes the familiar form. Using a weighted Witten identity, they proved a weighted positive mass theorem for spin manifolds~\cite{BaldaufOzuch22}, and subsequently Chu and Zhu proved the weighted positive mass theorem in the non-spin case for $3\le n\le 7$ in the spirit of the Schoen--Yau minimal hypersurfaces argument~\cite{ChuZhu24}.

\begin{Theorem}[Baldauf--Ozuch \cite{BaldaufOzuch22}, Chu--Zhu \cite{ChuZhu24}]\label{thm-BO-CZ}
	Let $(M^n,g,{\rm e}^{-f}{\rm d}V_g)$ be an asymptotically flat weighted manifold with $R_f\in L^1(M)$. Assume either $3\le n\le 7$ or $M$ is spin.
Then if $R_f\ge 0$, $\mathfrak{m}_f(g)\ge 0$ with equality if and only if $(M^n,g)$ is isometric to Euclidean space and $f\equiv0$.
\end{Theorem}
The reader is directed to Definition~\ref{def-AF} for the assumed decay rates for asymptotically flat weighted manifolds.

Somewhat surprisingly, Law, Lopez, and Santiago \cite{LawLopezSantiago24} showed that after a particular conformal transformation, the weighted positive mass theorem reduces precisely to the usual (unweighted) positive mass theorem, which can then be applied directly. In fact, the result they obtain is actually stronger than the previously established more technical proofs.
\begin{Theorem}[Law--Lopez--Santiago \cite{LawLopezSantiago24}]\label{thm-LLS}
Let $(M^n,g,{\rm e}^{-f}{\rm d}V_g)$ be an asymptotically flat weight\-ed manifold with $f=o_2\bigl(|x|^{-(n-2)/2}\bigr)$ and $R_f\in L^1(M)$. Set
	\[
	\widetilde g = {\rm e}^{-\frac{2}{n-1}f} g.
	\]
	Then the scalar curvature of $\widetilde g$, $\widetilde R$ is given by
	\[
	\widetilde R = {\rm e}^{\frac{2}{n-1}f}\left(R_f + \frac{1}{n-1} |\nabla f|^2\right),
	\]
	and
	\[
	\mathfrak{m}_f(g) = \mathfrak{m}(\widetilde g).
	\]
	In particular, if $R_f\ge -\frac{1}{n-1}|\nabla f|^2$ then $R_{\widetilde g}\ge 0$ and hence $\mathfrak{m}_f(g)\ge 0$.
\end{Theorem}

Without this insight, one may be tempted to pursue a proof of a weighted Riemannian Penrose inequality by adapting the proofs of Huisken and Ilmanen \cite{HuiskenIlmanen01} or Bray~\cite{Bray01} to the weighted setting.\footnote{The present author at least was, prior to learning of the article~\cite{LawLopezSantiago24}.} However, it turns out that the expected weighted Riemannian Penrose inequality follows from the same conformal transformation then applying the standard Riemannian Penrose inequality to the conformal metric. In this note, we show this weighted Riemannian Penrose inequality and some related results on weighted manifolds that follow from this same conformal picture. In order to better understand why this approach works so well, we investigate the quantity
\[
S_f=R_f+\frac{1}{n-1}|\nabla f|^2.
\] In particular, we study the geometric mass invariant associated with $S_f$ via Michel's mass invariant formalism \cite{Michel11}, and establish that it is precisely the weighted mass.
\begin{Theorem}\label{thm-Michel}
	Let $(M^n,g,{\rm e}^{-f}{\rm d}V_g)$ be an asymptotically flat weighted manifold and define the densitised curvature operator $\Phi(g,f)=S_f(g) {\rm e}^{-f}{\rm d}V_g$. The kernel of $D\Phi^*_{(\mathring g,0)}$ consists of affine functions $V$ defining geometric mass invariants in the sense of Michel~{\rm \cite{Michel11}}
	\begin{equation*}
		\mathfrak{m}(V,g,f)=\lim_{\rho\to\infty}\int_{S_\rho}\mathbb U_i(V,g-\mathring g,f) \nu^i{\rm d}S,
	\end{equation*}
	where $\mathbb U(V,g-\mathring g,f)$ satisfies
	\[
		\mathrm{div}_{\mathring g}\mathbb U(V,g-\mathring g,f)=V  D\Phi_{(\mathring g,0)}(g-\mathring g,f).
	\]
	Furthermore, for $V=1$ and assuming $S_f\in L^1$ then we have
	\[
		\mathfrak{m}(1,g,f)=\mathfrak{m}_f(g).
	\]
\end{Theorem}

With this in mind, it is not surprising that we can establish so many properties of the weighted mass via the conformal metric $\widetilde g$, whose scalar curvature is~$S_f$. This framework also naturally leads to a definition of centre of mass for weighted manifolds, which agrees with the standard (unweighted) centre of mass of~$\widetilde g$.

Examining the linearisation of $S_f$ and its formal adjoint leads to a natural definition of weighted static metrics, which can be understood as candidates for minimisers of the weighted mass. We prove that this definition of weighted staticity can also be understood via the conformal metric $\widetilde g$. From this we prove the following weighted static uniqueness theorem for what we call the $f$-Schwarzschild family of metrics, which differ from regular Schwarzschild metrics by the conformal factor \smash{${\rm e}^{-\frac{2f}{n-1}}$} (see Definition~\ref{def-f-schwarzschild}).
\begin{Theorem}\label{thm-fstatic-bh-rigidity-intro}
	Let $(M^n,g,{\rm e}^{-f}{\rm d}V_g)$ be an asymptotically flat weighted manifold with compact boundary
	$\Sigma=\partial M$, and assume $f\to 0$ at infinity.
	Assume either $3\le n\le 7$ or $M$ is spin, and for $n>3$ further assume $\Sigma$ is connected.
Suppose there exists an $f$-static potential $V\in C^\infty(M)$ satisfying
	\[
	V>0\ \text{on }M\setminus\Sigma,\qquad V=0\ \text{on }\Sigma,\qquad V\to 1\ \text{at infinity}.
	\]
Then $(M^n,g,{\rm e}^{-f}{\rm d}V_g)$ is $f$-Schwarzschild of mass $m=\mathfrak{m}_f(g)$.
\end{Theorem}

Theorem~\ref{thm-fstatic-bh-rigidity-intro} is the weighted version of the classical static black hole uniqueness theorems, originally due to Israel \cite{israel1967event}, then improved by Bunting and Masood-ul-Alam \cite{BuntingMasoodUlAlam87} to allow for a~disconnected horizon and extended to higher dimensions by Gibbons, Ida and Shiromizu \cite{GibbonsIdaShiromizu02}.

Many other geometric quantities have natural weighted versions that fit together within a~coherent weighted manifold framework (see \cite[Table~1]{BaldaufOzuch22}, for example). This leads to natural definitions of weighted minimal surfaces, weighted outer-minimising surfaces (see Definition~\ref{def-weighted-min}), and to conjecture a weighted Riemannian Penrose inequality. However, we show this follows immediately from the conformal transformation discussed above.
\begin{Theorem}\label{thm-weightedPenrose-intro}
	Let $(M^n,g,{\rm e}^{-f}{\rm d}V_g)$ be a weighted asymptotically flat manifold of dimension $3\le n\le 7$ containing a closed $f$-outer-minimising $f$-minimal hypersurface $\Sigma$. Assume $S_f\ge 0$, then
	\begin{equation}\label{eq-weightedPenrose}
		\mathfrak{m}_f(g) \ge \frac12\left(\frac{A_f(\Sigma)}{\omega_{n-1}}\right)^{\frac{n-2}{n-1}}.
	\end{equation}
	Furthermore, equality holds in \eqref{eq-weightedPenrose} if and only if $(M^n,g,{\rm e}^{-f}{\rm d}V_g)$ is an $f$-Schwarzschild manifold.
\end{Theorem}

We remark that it appears as though much of what one might want to study in connection to this weighted mass ultimately can be reduced to the unweighted situation by means of this straightforward conformal transformation.

\section{The weighted mass}\label{sec-michel}

Throughout this note, $(M^n,g)$ denotes a smooth Riemannian manifold of dimension $n\geq3$. We say that $(M^n,g)$ is \emph{asymptotically flat} of order $\tau>0$ if there exists a compact set $K\subset M$, a~radius $R>0$, and a diffeomorphism
\[
\Psi \colon \  M\setminus K \longrightarrow \mathbb{R}^n\setminus B_R(0),
\]
on $M\setminus K$ such that
\[
g_{ij} - \mathring g_{ij} = O_2(r^{-\tau}),
\]
where $r=|x|$ is the standard radial coordinate in the end and $\mathring g$ is a background metric isometric to the Euclidean metric on $M\setminus K$.

\begin{Definition}\label{def-AF}
	A \emph{weighted manifold} (or manifold with density) is a triple $(M^n,g,{\rm e}^{-f}{\rm d}V_g)$, where $f\in C^\infty(M)$ and the reference measure is ${\rm e}^{-f}{\rm d}V_g$. We say a weighted manifold is asymptotically flat of order $\tau>\frac{n-2}{2}$ if $(M^n,g)$ is asymptotically flat of order $\tau$ and $f=O_2(r^{-\tau})$.
\end{Definition}

We remark that the weighted area of a hypersurface $A_f(\Sigma)=\int_\Sigma {\rm e}^{-f}{\rm d}\sigma$ has first variation $\int_\Sigma H_f {\rm e}^{-f}{\rm d}\sigma$ under normal variations, where $H_f=H-\partial_\nu f$ is the weighted mean curvature (cf.~\cite{BaldaufOzuch22}). This leads us to make the following natural definitions.
\begin{Definition} \label{def-weighted-min}
	We say a hypersurface in a weighted manifold $(M^n,g,{\rm e}^{-f}{\rm d}V_g)$ is \emph{$f$-minimal} if the weighted mean curvature $H_f=H-\partial_\nu f$ of the hypersurface vanishes. We further say that a~closed hypersurface $\Sigma$ in a weighted manifold is \emph{$f$-outer-minimising} if it minimises the weighted area over all hypersurfaces enclosing $\Sigma$.
\end{Definition}

With the preliminaries out of the way, we now turn to prove Theorem~\ref{thm-Michel}, defining the weighted mass as a natural mass invariant associated to a curvature map.

Let $(M^n,g,{\rm e}^{-f}{\rm d}V_g)$ be a weighted $n$-dimensional asymptotically flat manifold with decay $\tau>\frac{n-2}{2}$. We define a densitised curvature map
\[
\Phi(g,f)= S_f {\rm e}^{-f}{\rm d}V_g,
\]
where
\[
	S_f=R_f+\frac{1}{n-1}|\nabla f|^2=R+2\Delta(f)- \frac{n-2}{n-1}|\nabla f|^2,
\]
motivated by Theorem~\ref{thm-LLS}. Following Michel's formalism on geometric mass invariants \cite{Michel11}, there is a $1$-form valued differential operator $\mathbb U$ defined by
\begin{equation*}
	V D\Phi_{(\mathring{g},\mathring{f})}(h, \varphi)-\langle(h,\varphi),D\Phi_{(\mathring{g},\mathring{f})}^*(V)\rangle=\mathrm{div}_{\mathring g}\mathbb U(V,h,\varphi){\rm d}V_{\mathring g},
\end{equation*}
where $h=g-\mathring g$, $\varphi=f-\mathring f$, $V$ is a scalar function, and $\langle \cdot,\cdot\rangle$ is the natural pairing. We note here that instead of using the densitised map $\Phi$, one could instead define the adjoint map with respect to a weighted $L^2$ inner product. Although we include the background weight function $\mathring f$ in the expression above for clarity, we will be taking it to vanish identically in the definition of the weighted mass. Choosing $V$ in the kernel of \smash{$D\Phi_{(\mathring{g},\mathring{f}\equiv 0)}^*$} and integrating $\mathrm{div}_{\mathring g}\mathbb U$ over $M$ gives a flux integral at infinity, which corresponds to a geometric invariant \`a la Michel~\cite{Michel11}. Although for the mass definition we need only compute the variation around the base point $(\mathring g,0)$, the linearisation and its adjoint at a general point will be useful in what follows. We compute
\begin{align*}
		D\Phi_{(g,f)}(h,\varphi) ={}& \left(
		\bigl( \nabla^i \nabla^j h_{ij} - \Delta(h^i{}_i) \bigr)
		- 2\left( \nabla^j h_{ij} - \frac{1}{2}\nabla_i(h^k{}_k) \right)\nabla^i f\right. \nonumber\\
		&  - h_{ij}\left( \mathrm{Ric}^{ij} + 2\nabla^i \nabla^j f - \frac{n-2}{n-1} \nabla^i f \nabla^j f \right)\nonumber \\
		& \left. +\left(2 \Delta \varphi - \frac{2(n-2)}{n-1} \nabla^i f  \nabla_i \varphi\right)
		 +S_f(g)\left( \frac12 \operatorname{tr}_g(h)-\varphi \right)\right){\rm e}^{-f}{\rm d}V_g,
\end{align*}
with formal adjoint
\begin{gather}
		D_g\Phi_{(g,f)}^*(V)^{ij}
		  = \left(
		\nabla^i\nabla^j V-(\Delta V) g^{ij}
		+\langle\nabla V,\nabla f\rangle g^{ij}\vphantom{\frac{1}{n-1}}\right. \nonumber\\
\left. \hphantom{D_g\Phi_{(g,f)}^*(V)^{ij}=}{}
		- V\left(\mathrm{Ric}^{ij}+\nabla^i\nabla^j f+\frac{1}{n-1} (\nabla^i f)(\nabla^j f)\right)
 +\frac12 V S_f g_{ij}
		\right) {\rm e}^{-f}{\rm d}V_g,\label{eq-adjoint}\\
		D_f\Phi_{(g,f)}^*(V)
		  = \left(
		2\Delta V
		-\frac{2n}{n-1} \nabla_i V \nabla^i f
		-\frac{2}{n-1} V \Delta f
		 +\frac{2}{n-1} V |\nabla f|^2
		- V S_f
		\right) {\rm e}^{-f}{\rm d}V_g,\nonumber
\end{gather}
where indices are raised and lowered with $g$.

From this we see that linearising around $(\mathring g, \mathring f\equiv0)$ gives
\begin{equation*}
	D\Phi_{(\mathring g,0)}(h,\varphi)=\bigl( 	\mathring{\nabla}^i \mathring{\nabla}^j h_{ij} - \mathring{\Delta}(\operatorname{tr}_{\mathring g}(h)) +2\mathring{\Delta}(\varphi)\bigr){\rm d}V_{\mathring g}
\end{equation*}
and
\begin{gather*}
		D_g\Phi_{(\mathring g,0)}^*(V)^{ij}=   \bigl(\mathring{\nabla}^i \mathring{\nabla}^j V - (\mathring{\Delta} V)  \mathring g^{ij}\bigr){\rm d}V_{\mathring g},\\
		D_f\Phi_{(\mathring g,0)}^*(V)=   2 \mathring{\Delta} V {\rm d}V_{\mathring g}.
\end{gather*}
We see immediately that the kernel of $D\Phi^*_{(\mathring g,0)}$ consists of affine functions. Taking $V\equiv 1$, and using $h=g-\mathring g$ and $\varphi=f$, we have
\[
	\mathbb U(1,g-\mathring g,f)=\mathrm{div}_{\mathring g}(g)-\mathring{\nabla}\operatorname{tr}_{\mathring g}(g)+2\mathring{\nabla} f.
\]
We then immediately see that integrating this over large spheres recovers the weighted mass
\begin{equation*}
	\lim_{\rho\to\infty}\int_{S_\rho}\mathbb U_i(1,g-\mathring g,f){\rm d}S^i=2(n-1)\omega_{n-1}\mathfrak{m}_f(g),
\end{equation*}
making use of the fact that ${\rm e}^{-f}=1+O(r^{-\tau})$. That is, the weighted mass is the geometric mass invariant associated with $S_f{\rm e}^{-f}{\rm d}V_g$ and we have established Theorem~\ref{thm-Michel}. An important aspect of this formalism for mass invariants is that it gives a somewhat universal framework to prove that they are indeed geometric quantities, independent of the choice of coordinates at infinity~\cite{Michel11}. However, for the weighted mass this readily follows from the coordinate invariance of the ADM mass~\cite{BaldaufOzuch22,LawLopezSantiago24}, so in the interest of brevity we do not explicitly present this here.

We remark that one could insert other affine functions into the construction above and analogously to the regular ADM mass, define a weighted centre of mass. Namely, one might define the centre of mass vector of an asymptotically flat weighted manifold as
\begin{align*}
		\mathfrak c^a_f(g)={}&\frac{1}{\mathfrak{m}_f(g)}\lim_{\rho\to\infty}\int_{S_\rho}\mathbb U_i(x^a,g-\mathring g,f){\rm d}S^i\\		={}&\frac{1}{\mathfrak{m}_f(g)}\lim_{\rho\to\infty}\int_{S_\rho}x^a\bigl(\mathring{\nabla}^{j}g_{ij}-\mathring{\nabla}_i\big(\operatorname{tr}_{\mathring g}g\big)+2 \mathring{\nabla}_if\bigr)
		\\ &  -(g-\mathring g)_{ia}
		+\bigl(\operatorname{tr}_{\mathring g}(g-\mathring g)-2f\bigr)\delta_i^a {\rm d}S^i,
\end{align*}
where $x^a$ is a coordinate function. We next show that this is simply the same as the standard centre of mass for the conformal metric \smash{$\widetilde g = {\rm e}^{-\frac{2}{n-1}f} g$}.
\begin{Proposition}
	Let $(M^n,g,{\rm e}^{-f}{\rm d}V_g)$ be a smooth asymptotically flat weighted manifold of order $\tau>\frac{n-2}{2}$, and let \smash{$\widetilde g = {\rm e}^{-\frac{2}{n-1}f} g$}. Then
	\[
	\mathfrak c^a_f(g)=\mathfrak c^a_{\mathrm{ADM}}(\widetilde g ),
	\]
	where $\mathfrak c_{\rm ADM}$ is the usual centre of mass for an asymptotically flat manifold given by
	\begin{gather} \label{eq-CoM}
			\mathfrak c_{\rm ADM}^a(\widetilde g) =\lim_{\rho\to\infty}\frac{1}{\mathfrak{m}_{\rm ADM}(\widetilde g)}\int_{S_\rho}
			x^a\bigl(\mathring\nabla^j(\widetilde g)_{ij}-\mathring\nabla_i\operatorname{tr}_{\mathring g}(\widetilde g)\bigr)
			-(\widetilde g -\mathring g)_{ia}+\operatorname{tr}_{\mathring g}(\widetilde g-\mathring g) \delta_i^{ a}{\rm d}S^i,
	\end{gather}
	whenever $\mathfrak c_{\rm ADM}(\widetilde g)$ is well defined.\footnote{ To ensure the centre of mass is well defined requires an additional parity assumption \cite{ReggeTeitelboim74}.}
\end{Proposition}
\begin{proof}
	The decay conditions allow us to write
	\begin{equation*}
		\widetilde g
		= g - \frac{2}{n-1}f \mathring g + O_2\bigl(r^{-2\tau}\bigr).
	\end{equation*}
	Inserting this into \eqref{eq-CoM} and recalling $\mathfrak{m}_{\rm ADM}(\widetilde g)=\mathfrak{m}_f(g)$ one quickly arrives at
	\[
	\mathfrak c^a_f(g)=\mathfrak c^a_{\mathrm{ADM}}(\widetilde g). \tag*{\qed}
	\]\renewcommand{\qed}{}
\end{proof}

Although this weighted centre of mass does turn out to be simply the regular centre of mass of a conformal metric, we include it partially to complete the mass invariant picture but also for the following reason. The usual centre of mass of an asymptotically flat manifold can be instead defined in terms of a foliation of the asymptotic end by constant mean curvature (CMC) surfaces \cite{HuiskenYau96}. However, pulling back these CMC surfaces via the conformal transformation used above does not give a foliation of constant weighted mean curvature $H_f$ surfaces. It seems reasonable though to believe by analogy that whenever the weighted mass is positive, we can find a foliation of the end of a weighted manifold by surfaces of constant $H_f$ and centered at $\mathfrak c_f$. However, we do not pursue this question further here.

\section{Weighted staticity}\label{sec-f-static}

The preceding section defined the weighted mass as the geometric invariant at infinity associated to the map $\Phi(g,f)=S_f{\rm e}^{-f}{\rm d}V_g$. By direct analogy to the usual notion of static metrics, we make the following definition.

\begin{Definition}[$f$-static potential]\label{def-f-static}
	A nontrivial function $V\in C^\infty(M)$ is called a \emph{weighted static potential} or an \emph{$f$-static potential} if it satisfies
	\[
	D\Phi_{(g,f)}^*(V)=0.
	\]
	We say $(g,f)$ is \emph{$f$-static} if it admits an $f$-static potential.
\end{Definition}
For simplicity, we make use of the de-densitised adjoint map $F^*$ given by
\begin{equation*}
	D_g\Phi^*_{(g,f)}(V)=F_g^*(V) {\rm e}^{-f}{\rm d}V_g,
	\qquad
	D_f\Phi^*_{(g,f)}(V)=F_f^*(V) {\rm e}^{-f}{\rm d}V_g.
\end{equation*}
We use the same conformal change as earlier, and consider the metric
\[
\widetilde g = {\rm e}^{-\frac{2}{n-1}f} g.
\]
In this section, we show that if $g$ is $f$-static with $f$-static potential $V$, then $\widetilde g$ is static with static potential
\[
u={\rm e}^{-\frac{f}{n-1}}V.
\]
We begin by calculating some identities for $u$ and the conformal metric $\widetilde g$.

\begin{Lemma}\label{lem-fstatic-conf-identities}
	With $\widetilde g$, $u$ and $F$ as above, we have
	\begin{align}
		\label{eq-conf-lap}
&		\widetilde\Delta u
={\rm e}^{-\frac{f}{n-1}}\left(\frac12 F_f^*(V)+\frac12 V S_f\right),\\
		\label{eq-conf-hessric}
&		\widetilde\nabla^2 u-u \widetilde{\mathrm{Ric}}
={\rm e}^{\frac{f}{n-1}}\left(F_g^*(V)+\frac12 F_f^*(V) g\right).
	\end{align}
\end{Lemma}

\begin{proof}
	Writing \smash{$\psi=-\frac{f}{n-1}$} for the sake of exposition, we have the standard formulas for a~conformal metric \smash{$\widetilde g= {\rm e}^{2\psi}g$},
	\begin{align}
		&\widetilde\Delta w
={\rm e}^{-2\psi}\bigl(\Delta w+(n-2) \nabla_i\psi \nabla^i w\bigr), \label{eq-conf-proof-lap}\\
&		\widetilde\nabla^2_{ij} w
=\nabla_{ij}^2 w-\nabla_iw \nabla_j\psi-\nabla_i\psi \nabla_jw+\nabla_k\psi \nabla^k w g_{ij}, \label{eq-conf-proof-hess}\\
& \widetilde{\mathrm{Ric}}_{ij}
=\mathrm{Ric}_{ij}-(n-2)\big(\nabla_{ij}^2\psi-\nabla_i\psi \nabla_j\psi\big)
-\big(\Delta\psi+(n-2)|\nabla\psi|^2\big)g_{ij}. \label{eq-conf-proof-ric}
	\end{align}
	The formula \eqref{eq-conf-lap} can be derived from \eqref{eq-conf-proof-lap} as follows. Using $u={\rm e}^{\psi}V$, we have
	\[
	\nabla u={\rm e}^{\psi}(\nabla V+V\nabla\psi),
	\]
	and then
	\[
	\Delta u
	={\rm e}^{\psi}\bigl(\Delta V+2\nabla_i\psi \nabla^i V+V\Delta\psi+V|\nabla\psi|^2\bigr).
	\]
	Moreover,
	\[
	\nabla_i\psi \nabla^i u
	={\rm e}^{\psi}\bigl(\nabla_i\psi \nabla^i V+V|\nabla\psi|^2\bigr).
	\]
	Inserting these into \eqref{eq-conf-proof-lap} gives
	\begin{equation}\label{eq-conf-proof-lap-expanded-new}
		\widetilde\Delta u
		={\rm e}^{-\psi}\bigl(\Delta V+n\nabla_i\psi \nabla^i V+V\Delta\psi+(n-1)V|\nabla\psi|^2\bigr).
	\end{equation}
	We next show that the terms in parentheses can be expressed in terms of $F_f^*(V)$ and $S_f$. In particular, from \eqref{eq-adjoint} we have
	\[
	\frac12 F_f^*(V)
	=\Delta V+n\nabla_i\psi \nabla^i V+V\Delta\psi+(n-1)V|\nabla\psi|^2-\frac12 V S_f,
	\]
	which is finally substituted into \eqref{eq-conf-proof-lap-expanded-new} to give \eqref{eq-conf-lap}.
	
	We next turn to establish \eqref{eq-conf-hessric} using \eqref{eq-conf-proof-hess} and \eqref{eq-conf-proof-ric}. First note
	\[
	\nabla_{ij}^2 u
	={\rm e}^{\psi}\bigl(\nabla_{ij}^2 V + \nabla_iV \nabla_j\psi + \nabla_i\psi \nabla_jV + V\nabla_{ij}^2\psi + V \nabla_i\psi \nabla_j\psi\bigr),
	\]
	so using \eqref{eq-conf-proof-hess} we have
	\begin{align}
		\widetilde\nabla_{ij}^2 u
		&={\rm e}^{\psi}\bigl(\nabla_{ij}^2 V + V\nabla_{ij}^2\psi - V \nabla_i\psi\nabla_j \psi
		+\big(\nabla_k V \nabla^k\psi+V|\nabla\psi|^2\big)g_{ij}\bigr). \label{eq-conf-proof-hess-u}
	\end{align}
	Multiplying \eqref{eq-conf-proof-ric} by $u$ and subtracting it from \eqref{eq-conf-proof-hess-u} gives
	\begin{align*}
		\widetilde\nabla_{ij}^2 u-u \widetilde{\mathrm{Ric}}_{ij}
={}&{\rm e}^{\psi}\bigl(\nabla_{ij}^2 V - V\mathrm{Ric}_{ij}
		+(n-1)V(\nabla_{ij}^2\psi-\nabla_i\psi \nabla_j\psi) \nonumber\\
		&
		+\big(\nabla_k V \nabla^k\psi+V\Delta\psi+(n-1)V|\nabla\psi|^2\big)g_{ij}\bigr). \label{eq-conf-proof-hessric-new}
	\end{align*}
	Finally, writing everything in terms of $f$ instead of $\psi$ then comparing with \eqref{eq-adjoint} shows that the term in parentheses is exactly
	\[
	F_g^*(V)+\frac12 F_f^*(V) g. \tag*{\qed}
	\]
\renewcommand{\qed}{}
\end{proof}

\begin{Proposition}\label{prop-fstatic-to-static}
	Suppose $V$ is an $f$-static potential for $(M^n,g,{\rm e}^{-f}{\rm d}V_g)$ then $(\widetilde g,u)$ defined above is static vacuum.
\end{Proposition}

\begin{proof}
	Notice that by Lemma \ref{lem-fstatic-conf-identities}, the condition that $\widetilde g$ be static with static potential $u$ is equivalent to $S_f\equiv0$ and $g$ is $f$-static with weighted static potential $V$. So we only must show that if $g$ is $f$-static then $S_f$ must identically vanish. To this end, assume $F_g^*(V)=0$ and $F_f^*(V)=0$. Tracing $F_g^*(V)$ gives
	\[
	(1-n)\Delta(V)+n\nabla_k(V)\nabla^k(f)-V\left(R+\Delta(f)+\frac{1}{n-1}|\nabla f|^2\right)+\frac{n}{2}VS_f,
	\]
	which we compare to
	\[
	\frac{n-1}{2}F_f^*(V)=(n-1)\Delta V-n\nabla_k(V)\nabla^k(f)-V\Delta f+V|\nabla f|^2-\frac{n-1}{2}V S_f.
	\]
	Adding these together gives
	\begin{align*}
	&	0=\operatorname{tr}_g F_g^*(V)+\frac{n-1}{2}F_f^*(V)
		 =V\left(-R-2\Delta f+\frac{n-2}{n-1}|\nabla f|^2+\frac{1}{2}S_f \right)
		 =-\frac12 V S_f.
	\end{align*}
	It remains to show that $V$ cannot vanish on an open set, which follows by virtue of the fact that $V$ satisfies an elliptic equation. In fact, by \eqref{eq-conf-hessric}, \smash{$u={\rm e}^{-\frac{f}{n-1}}V$} satisfies $\widetilde \Delta u-\widetilde R u=0$ so by unique continuation, neither $u$ nor $V$ can vanish on an open set unless $V\equiv0$. It follows that $S_f\equiv0$.
\end{proof}

Since we can identify $f$-static metrics with static metrics via a conformal transformation, static uniqueness results translate directly over to this setting with the unique $f$-static metric being conformal to a Schwarzschild metric.

\begin{Definition}\label{def-f-schwarzschild}
	Fix $m>0$, and let $r_0=\big(\frac{m}{2}\big)^{\frac{1}{n-2}}$ and
	$M_m=\mathbb{R}^n\setminus B_{r_0}(0)$.
	Given any $f\in C^\infty(M_m)$ with $f\to 0$ as $r\to\infty$, define the metric
	\[
	g_{m,f}
	={\rm e}^{\frac{2}{n-1}f}\left(1+\frac{m}{2 r^{ n-2}}\right)^{\frac{4}{n-2}}\delta
	\quad\text{on }M_m,
	\]
	where $\delta$ denotes the standard Euclidean metric.
	
	We call $(M_m,g_{m,f},{\rm e}^{-f}{\rm d}V_g)$ an \emph{$f$-Schwarzschild manifold of weighted mass $m$}.
\end{Definition}

Theorem~\ref{thm-fstatic-bh-rigidity-intro} readily follows.

\begin{proof}[Proof of Theorem~\ref{thm-fstatic-bh-rigidity-intro}]
	Let $\widetilde g$ and $u$ be as above. By Proposition \ref{prop-fstatic-to-static}, $\widetilde g$ is static with static potential $u$ such that the zero sets of $u$ and $V$ are identical, $\{u=0\}=\{V=0\}=\Sigma$, and $u\to 1$ at infinity. Moreover $u>0$ on $M\setminus\Sigma$ and $u=0$ on $\Sigma$, so in the language of general relativity, $\Sigma$~is a non-degenerate horizon for $(\widetilde g,u)$. That is, $(M,\widetilde g,u)$ is an asymptotically flat vacuum static triple with compact horizon $\Sigma=\{u=0\}$ and
	$u\to 1$ at infinity.
	
	We therefore simply apply the classical black hole uniqueness results. For example, if $n=3$ the work of Bunting and Masood-ul-Alam \cite{BuntingMasoodUlAlam87} applies, and if $n>3$ but $\Sigma$ is connected, then the work of Gibbons, Ida, and Shiromizu \cite{GibbonsIdaShiromizu02} applies. In either case $\widetilde g$ is a positive mass Schwarzschild metric and therefore $(M^n,g,{\rm e}^{-f}{\rm d}V_g)$ is
	$f$-Schwarzschild of mass $m$.
\end{proof}

\section{The weighted Riemannian Penrose inequality}

In this section, we prove a weighted Riemannian Penrose inequality with the rigidity case corresponding to the $f$-Schwarzschild family of metrics.

\begin{Lemma}\label{lem-fmin-to-min}
	Let $(M^n,g,{\rm e}^{-f}{\rm d}V_g)$ be a smooth weighted manifold with $\Sigma\subset M$ a closed hypersurface, and set
	\[
	\widetilde g = {\rm e}^{-\frac{2}{n-1}f} g.
	\]
	Then,
	\begin{enumerate}\itemsep=0pt
		\item[$(1)$] $\Sigma$ is $f$-minimal in $(M^n,g,{\rm e}^{-f}{\rm d}V_g)$ if and only if $\Sigma$ is minimal in $(M,\widetilde g)$;
		\item[$(2)$] the $\widetilde g$-area equals the weighted area,
		\[
		A_{\widetilde g}(\Sigma)=\int_\Sigma {\rm d}A_{\widetilde g}=\int_\Sigma {\rm e}^{-f} {\rm d}A_g = A_f(\Sigma);
		\]
		\item[$(3)$] $\Sigma$ is $f$-outer-minimising in $(M^n,g,{\rm e}^{-f}{\rm d}V_g)$ if and only if $\Sigma$ is outer-minimising in $(M,\widetilde g)$.
	\end{enumerate}
\end{Lemma}

\begin{proof}
	Letting \smash{$\widetilde g={\rm e}^{-\frac{2f}{n-1}}g$}, the mean curvature transforms as
	\[
	\widetilde H ={\rm e}^{\frac{f}{n-1}}\big(H-\partial_\nu f\big)={\rm e}^{\frac{f}{n-1}}H_f,
	\]
	so (1) holds. The relationship between areas (2) follows even more straightforwardly from the conformal transformation formula.
	
	Finally, since we have $A_{\widetilde g}(\Sigma')=A_f(\Sigma')$ for any $\Sigma'$ enclosing $\Sigma$, and the choice of metric plays no role in whether $\Sigma'$ encloses $\Sigma$, (3) holds.
\end{proof}

We now prove Theorem~\ref{thm-weightedPenrose-intro}.

\begin{proof}[Proof of Theorem~\ref{thm-weightedPenrose-intro}]
	Let $\widetilde g={\rm e}^{-\frac{2}{n-1}f}g$ be the same conformal metric as earlier.
	By Theorem~\ref{thm-LLS}, the scalar curvature satisfies
	\[
	\widetilde R = {\rm e}^{\frac{2}{n-1}f}S_f\geq0,
	\qquad \text{and} \qquad
	\mathfrak{m}(\widetilde g)=\mathfrak{m}_f(g).
	\]
	
	Since $\Sigma$ is $f$-minimal and $f$-outer-minimising in $(M^n,g,{\rm e}^{-f} {\rm d}V_g)$, Lemma \ref{lem-fmin-to-min} implies that
	$\Sigma$ is minimal and outer-minimising in $(M,\widetilde g)$, and that $A_{\widetilde g}(\Sigma)=A_f(\Sigma)$.
	We may therefore apply the Riemannian Penrose inequality \cite{Bray01,BrayLee09, HuiskenIlmanen01} to $(M,\widetilde g)$ to obtain
	\[
	\mathfrak{m}(\widetilde g) \ge \frac12\biggl(\frac{A_{\widetilde g}(\Sigma)}{\omega_{n-1}}\biggr)^{\frac{n-2}{n-1}}
	=\frac12\biggl(\frac{A_f(\Sigma)}{\omega_{n-1}}\biggr)^{\frac{n-2}{n-1}}.
	\]
	Since $\mathfrak{m}(\widetilde g)=\mathfrak{m}_f(g)$, \eqref{eq-weightedPenrose} follows.
	
	Furthermore, rigidity in the Riemannian Penrose inequality\footnote{In \cite{BrayLee09}, spin is assumed for rigidity, but this restriction was removed by subsequent work of McFeron and Sz{\'e}kelyhidi \cite{mcferon2012positive}. I would like to thank an anonymous reviewer for pointing this out to me.} implies rigidity here. That is, $(M,\widetilde g)$ is isometric to a Schwarzschild exterior and therefore~$g$ is an $f$-Schwarzschild metric in the sense of Definition~\ref{def-f-schwarzschild}.
\end{proof}

We remark that in dimension $3$, there is in fact a stronger inequality than the Riemannian Penrose inequality. Namely the inequality between the Hawking mass and the ADM mass \cite{HuiskenIlmanen01}. It turns out we get an analogous result applying essentially the proof of Theorem~\ref{thm-weightedPenrose-intro} again.

Define a weighted Hawking mass (in dimension $3$) as
\begin{equation}\label{eq-weighted-hawking-formula}
	m_{H,f}(\Sigma;g)
	=\biggl(\frac{A_f(\Sigma)}{16\pi}\biggr)^{\frac12}
	\biggl(1-\frac{1}{16\pi}\int_\Sigma H_f^{ 2} {\rm d}A_g\biggr),
\end{equation}
by analogy to the usual Hawking mass. One can quickly check this is the same as the Hawking mass $\mathfrak{m}_H(\Sigma;\widetilde g)$ of $\Sigma$ in $(M,\widetilde g)$, and therefore applying the inequality
\begin{equation}\label{eq-Hawking-ADM}
	\mathfrak{m}_H(\Sigma;\widetilde g)\leq \mathfrak{m}_{\rm ADM}(\widetilde g)
\end{equation}
for outer-minimising surfaces \cite{HuiskenIlmanen01} and following the same proof as that of Theorem~\ref{thm-weightedPenrose-intro} gives a~weighted version of \eqref{eq-Hawking-ADM}.

\begin{Theorem}\label{thm-weighted-hawking-bound}
	Let $\bigl(M^3,g,{\rm e}^{-f}{\rm d}V_g\bigr)$ be asymptotically flat and assume
	\[
	S_f \ge 0.
	\]
	Let $\Sigma\subset M$ be a connected closed surface that is $f$-outer-minimising in $\bigl(M^3,g,{\rm e}^{-f}{\rm d}V_g\bigr)$.
	Then
	\begin{equation*}
		\mathfrak{m}_f(g) \ge \mathfrak{m}_{H,f}(\Sigma;g).
	\end{equation*}
\end{Theorem}
Note that the definition given by \eqref{eq-weighted-hawking-formula} does not use the weighted area form in the Willmore term. Perhaps a more natural definition of weighted Hawking mass would replace ${\rm d}A_g$ with ${\rm e}^{-f}{\rm d}A_g$, but it isn't clear that this definition relates to the weighted mass at infinity in the same way.

\subsection*{Acknowledgements}

The author is supported by Olle Engkvists stiftelse, the Magnus Bergvall Foundation, and foundations managed by The Royal Swedish Academy of Sciences. The author would also like to thank the anonymous referees for valuable comments that improved the clarity of this note.


\pdfbookmark[1]{References}{ref}
\LastPageEnding

\end{document}